# Soluciones débiles de las ecuaciones de Navier-Stokes mediante operadores acotados

Carlos de Jesús Prieto Chávez

Instituto Tecnológico de Celaya

4 de diciembre de 2013

# Índice

# 1. Introducción

En esta investigación se exponen las herramientas básicas para la demostración de lo que se denominan "Soluciones Turbulentas"de las ecuaciones de Navier-Stokes para un fluido viscoso incompresible.
La deducción de estas expresiones matemáticas son mostradas en varios libros de texto, sin embargo la demostración que se proporciona es enfocada al análisis tensorial e involucra a la Mecánica del Medio Continuo que no es más que un concepto matemático que define el campo de acción de todo fenómeno que posee deformaciones. La rigurosidad de las demostraciones están sustentadas mediante resultados preliminares de la Mecánica de Fluidos.

# 2. Planteamiento del Problema

El comportamiento de un fluido, de manera general, se puede clasificar como laminar y turbulento. La resolución de estas expresiones matemáticas responderían a la siguiente pregunta: *Para un movimiento laminar con condiciones iniciales dadas, ¿el comportamiento del fluido es laminar en todo instante?*
El planteamiento en el que se basa esta investigación es que si se conocen las condiciones iniciales para el movimiento de un fluido viscoso incompresible ¿Cuál será la respuesta de salida de este fluido en el dominio turbulento? ¿Será turbulento siempre?, por lo tanto es necesario saber como es está respuesta en el dominio del tiempo.

## 2.1. Objetivos

Se busca ampliar las herramientas matemáticas, las cuales permitan resolver problemas variados de la mecánica tanto en la teoría viscosa como en la mecánica del medio continuo. Se tratará de relacionar lo que es el análisis funcional con el tratado formal de la mecánica del medio continuo.

# 3. Antecedentes Históricos

La teoría de la viscosidad permite determinar el movimiento de un fluido mediante las ecuaciones de Navier-Stokes. Es necesario constatar que se puede modelar el movimiento de este fluido, estableciendo es siguiente *teorema de existencia:* hay solución a las ecuaciones de Navier-Stokes para un estado de velocidad arbitraria en un instante inicial. Esto fue lo que trato de probar Oseen [Ose10]. En su demostración,logró probar la existencia de esta solución por un corto tiempo después del instante dado de la condición inicial. También se puede verificar que la energía del fluido permanece acotada, pero no es posible deducir de este hecho que el movimiento permanece uniforme.En esta investigación se mostrará debido a qué, el movimiento se vuelve no uniforme en un tiempo finito.
De hecho, no es paradójico suponer que lo que regulariza el movimiento de disipación de energía no es suficiente para mantener las segundas derivadas de las componentes de velocidad acotadas y continuas. La teoría clásica de Navier y Stokes asume que las segundas derivadas son continuas y acotadas. Oseen sostuvo como hipótesis,en su trabajo sobre "Hidrodinámica", que no era un fenómeno natural. El mismo mostró como el movimiento obedecía las leyes de la mecánica, y las expreso en forma de ecuaciones integro-diferenciales, en las cuales solo contiene primeras derivadas de componentes de velocidad. Aquí muestro las ecuaciones de Oseen complementadas con una expresión(desigualdad) de disipación de energía. Estas expresiones pueden ser deducidas utilizando las ecuaciones de Navier-Stokes clásicas, y operando mediante integración por partes lo que proporciona que las derivadas de orden superior desaparezcan. Si no se logra demostrar el teorema de existencia propiamente dicho anteriormente, he probado lo siguiente: las relaciones en cuestión siempre tienen ."ªl menosüna solución para un estado inicial de velocidad dado, la cual está determinada para un tiempo de aplicación infinito en el que el origen es el estado inicial. En este trabajo se hace hincapié, en que al no poder acotar las derivadas de orden superior, son soluciones auto-satisfacidas por las ecuaciones de Navier-Stokes, por lo tanto son conocidas como soluciones turbulentas.

Una solución *turbulenta* tiene la siguiente estructura: se conforma de una sucesión de soluciones regulares.

Si al construir las soluciones para las ecuaciones de Navier-Stokes, se convierten en irregulares, entonces se puede asumir que existe soluciones que no son una composición de soluciones regulares. Del mismo modo, si la proposición es falsa,la semántica de *solución turbulenta* no pierde interés aunque por ahora no sea utilizada para describir los fluidos viscosos: sirve para representar problemas de la física-matemática para el cual la causa física de regularidad del sistema, no es suficiente para describir la regularidad total del propio sistema en todo el tiempo a la hora de escribir el modelo matemático.

En este trabajo se considera al fluido como un objeto de viscosidad ilimitada.

## 4. Consideraciones Preliminares

### 4.1. Notación

Usaremos la letra $\Pi'$ para un dominio arbitrario de puntos en el espacio. $\Pi'$ puede ser el espacio completo denotado por $\Pi$. $\bar{\omega}$ representa un dominio acotado en $\Pi$ el cual posee una frontera regular $\sigma$. Representaremos un punto arbitrario en el espacio $\Pi$ como $x$, con coordenadas euclideanas $x_i$ ($i = 1, 2, 3$) y distancia al origen $r_0$, y genera un elemento de volumen $\delta x$ y elementos de superficie $\delta x_1, \delta x_2, \delta x_3$. De la misma manera definimos a $y$ como un segundo punto arbitrario de $\Pi$. $r$ siempre denotará la distancia entre los puntos $x$ e $y$. Se utilizará el convenio de sumación de Einstein. En el presente trabajo se utilizo la letra $A$ para denotar constantes en el que no se especificó un valor.

Se usarán letras mayúsculas para representar funciones medibles, y para funciones continuas con primeras derivadas parciales continuas se ocupará letras minúsculas.

### 4.2. Desigualdad de Schwarz

$$\left[\iiint_{\Pi'} U(x)V(x)\delta x\right]^2 \leq \iiint_{\Pi'} U^2(x)\delta x \times \iiint_{\Pi'} V^2(x)\delta x \qquad (1)$$

El lado izquierdo de la desigualdad esta definido, si y solo si, el lado derecho es finito. Para esta sección es crucial esta desigualdad.

*Primera Aplicación:*

Si

$$U(x) = V_1(x) + V_2(x)$$

entonces

$$\sqrt{\iiint_{\Pi'} U^2(x)\delta x} \leq \iiint_{\Pi'} V_1^2(x)\delta x + \iiint_{\Pi'} V_2^2(x)\delta x$$

más general para $t$ constante

$$U(x) = \int_0^t V(x, t')dt'$$

entonces

$$\sqrt{\iiint_{\Pi'} U^2(x)\delta x} \leq \int_0^t dt' \sqrt{\iiint_{\Pi'} V^2(x, t0)\delta x} \qquad (2)$$

la parte izquierda de la desigualdad es finita cuando así lo sea la derecha.

*Segunda Aplicación:*

Consideremos $n$ constantes $\lambda_p$ y $n$ vectores constantes $\overrightarrow{\alpha_p}$. Se denota a $x + \overrightarrow{\alpha_p}$ como la traslación de $x$ debido a $\overrightarrow{\alpha_p}$. Tenemos

$$\iiint_\Pi \left[\sum_{p=1}^{p=n} \lambda_p U(x+\overrightarrow{\alpha_p})\right]^2 \delta x \leq \left[\sum_{p=1}^{p=n} |\lambda_p|\right]^2 \times \iiint_\Pi U^2(x)\delta x$$

es sencillo demostrar la expresión anterior, basta con desarrollar los términos cuadráticos y utilizando después la desigualdad de Schwarz. De la expresión, se puede deducir una relación útil. Sea $H(z)$ una función. Sí denotamos a $H(y-x)$ una función que se obtiene sustituyendo para coordenadas $z_i$ de $z$ los componentes $y_i - x_i$ del vector $\overrightarrow{xy}$, tenemos entonces

$$\iiint_\Pi \left[\iiint_\Pi H(y-x)U(y)\delta y\right]^2 \leq \left[\iiint_\Pi |H(z)|\right]^2 \times \iiint_\Pi U^2(y)\delta y; \tag{3}$$

la parte izquierda es finita si y solo si, las dos integrales de la parte derecha lo son.

### 4.3. Convergencia Fuerte en la Media ([Rie10])

Definición[1]: Se dice que una infinidad de funciones $U^*(x)$ posee una función $U(x)$ como límite fuerte de la media en el dominio $\Pi'$, cuando:

$$lim \iiint_{\Pi'} [U^*(x) - U(x)]^2 \delta x = 0. \tag{4}$$

Para cualquier función cuadrado sumable $A(x)$ en $\Pi'$ se tiene

$$lim \iiint_{\Pi'} U^*(x)A(x)\delta x = \iiint_{\Pi'} U(x)A(x)\delta x \tag{5}$$

De (4) y (5)

$$lim \iiint_{\Pi'} U^{*2}(x)\delta x = \iiint_{\Pi'} U^2(x)\delta x \tag{6}$$

*Convergencia Débil en la Media*
Definición: Un conjunto de infinitas funciones $U^*(x)$, tiene una función $U(x)$ como límite de convergencia en la media, cuando las siguientes dos proposiciones se cumplen:
a) el conjunto de números $\iiint_{\Pi'} U^{*2}(x)\delta x$ es acotado
b) para toda función cuadrado sumable $A(x)$ en $\Pi'$

$$lim \iiint_{\Pi'} U^*(x)A(x)\delta x = \iiint_{\Pi'} U(x)A(x)\delta x$$

*Ejemplo 1:* la secuencia $\sin x_1, \sin 2x_1, \sin 3x_1, \ldots$ converge débilmente a cero en todo el dominio $\bar{\omega}$
*Ejemplo 2:* Sea una infinidad de funciones $U^*(x)$ en un dominio $\Pi'$ converge al menos en una función $U(X)$. Esa función es límite en la media, cuando la cantidad de números $\iiint_{\Pi'} U^{*2}(x)\delta x$ sean acotados.

$$lim \iiint_{\Pi'_1} \iiint_{\Pi'_2} A(x)U^*(x)V^*(y)\delta x\delta y = \iiint_{\Pi'_1} \iiint_{\Pi'_2} A(x)U(x)V(y)\delta x\delta y \tag{7}$$

Donde $U^*(x)$ converge débilmente a $U(x)$ en el dominio $\Pi'_1$, y $V^*(y)$ converge débilmente a $V(x)$ y la integral

$$\iiint_{\Pi'_1} \iiint_{\Pi'_2} A^2(x)\delta x\delta y$$

---
[1]En la literatura actual el decir en la media está en desuso, sin embargo se utilizará para este artículo

es finita. Se tiene también

$$lim \iiint_{\Pi'} A(x)U^*(x)V^*(x)\delta x = \iiint_{\Pi'} A(x)U(x)V(x)\delta x \qquad (8)$$

donde en $\Pi'$, $A(x)$ es acotada, donde $U(x)$ es un límite fuerte de $U^*(x)$, y $V(x)$ es un límite débil de $U^*(x)$.
Es evidente que, si las funciones $U^*(x)$ convergen débilmente en la media a $U(x)$ en el dominio $\Pi'_1$

$$lim \left\{ \iiint_{\Pi'} [U^*(x) - U(x)]^2 \delta x - \iiint_{\Pi'} U^{*2}(x)\delta x + \iiint_{\Pi'} U^2(x)\delta x \right\} = 0$$

de la cual se obtiene la siguiente desigualdad

$$\iiint_{\Pi'} U^2(x)\delta x \leq lim\ inf \iiint_{\Pi'} U^{*2}(x)\delta x \qquad (9)$$

y el criterio para *convergencia fuerte:*
Las funciones $U^*(x)$ convergen fuertemente en la media del dominio $\Pi'$, cuando:

$$lim\ sup \iiint_{\Pi'} U^{*2}(x)\delta x \leq \iiint_{\Pi'} U^2(x)\delta x \qquad (10)$$

*Teorema de F. Riesz:* Una cantidad infinita de funciones $U^*(x)$ posee un límite débil en la media sobre el dominio $\Pi'$ si las siguientes condiciones se satisfacen:
a) Si el conjunto de números $\iiint_{\Pi'} U^{*2}(x)\delta x$ es acotado.
b) Para cada función cuadrado sumable $A(x)$ en $\Pi'$, las cantidades $\iiint_{\Pi'} U^*(x)A(x)\delta x$ poseen un único valor límite.
La condición b) puede ser remplazada como sigue:
b') Para cada cubo $c$ con lados paralelos a los ejes coordenados y vertices rotacionales, las catidades $\iiint_c U^*(x)\delta x$ poseen un único límite.
La demostración de este teorema, se realiza mediante la teoría de Lebesgue.

### 4.4. Método de diagonalización de Cantor

Consideremos un conjunto contable de cantidades cada una dependiente de un índice entero $n$: $a_n, b_n, \ldots (n = 1, 2, 3, \ldots)$. Supongamos que $a_n$ está acotado, $b_n$ está acotado, etc. El método de diagonalización de Cantor nos permite encontrar una secuencia de números $m_1, m_2, \ldots$, tal que cada una de las secuencias $a_{m_1}, a_{m_2}, \ldots$; $b_{m_1}, b_{m_2}, \ldots; \ldots$, convergen a un límite.
Aplicación: El siguiente resultado es consecuencia de los párrafos anteriores.
*Teorema fundamental de F. Riesz:* Sea una cantidad infinita de funciones $U^*(x)$ sobre un dominio $\Pi'$ tal que las cantidades $\iiint_{\Pi'} U^{*2}(x)\delta x$. Entonces siempre se puede extraer una secuencia la cual posee un límite débil en la media.
La condición a) se satisface, y el método de diagonalización de Cantor nos permite construir esa secuencia de funciones $U^*(x)$, la cual satisface la condición b').

### 4.5. Modos de continuidad de una función respecto a un parámetro

Sea una función $U(x, t)$ dependiente de un parámetro $t$. Se dice que esta función es *uniformemente continua* cuando las siguientes condiciones se cumplen:
a) si es continua con respecto a $x_1, x_2, x_3, t$
b) para cada valor particular $t_0$ de $t$, el máximo de $U(x, t_0)$ es finito

c) dado un número positivo $\epsilon$, se puede encontrar un número positivo $\eta$, tal que la desigualdad $|t-t_0| < \eta$ implica:

$$|U(x,t) - U(x,t_0)| < \epsilon$$

El máximo de $|U(x,t)|$ en $\Pi$, es entonces una función continua de $t$. Se dice que $U(x,t)$ es *fuertemente continua* en $t$ cuando, para cada valor particular $t_0$ de $t$, $\iiint_\Pi U^*(x,t_0)\delta x$ es finito, y para cada $\epsilon$ existe un $\eta$ tal que la desigualdad $|t - t_0| < \eta$ implica que:

$$\iiint_\Pi [U(x,t) - U(x,t_0)]^2 \delta x < \epsilon$$

La integral $\iiint_\Pi U^2(x,t)\delta x$ es entonces una función continua de $t$.

## 4.6. Relaciones entre una función y sus derivadas

Consideremos dos funciones $a(x)$ y $b(x)$ con primeras derivadas continuas. Sea $s$ la superficie de una esfera $S$ con centro en el origen donde el radio $r_0$ puede ser arbitrariamente largo. Sea

$$\varphi(r_0) = \iint_s u(x)a(x)\delta x_i$$

Tenemos que

$$\varphi(r_0) = \iiint_S \left[u(y)\frac{\partial a(y)}{\partial y_i} + \frac{\partial u(y)}{\partial y_i}a(y)\right]\delta y$$

La segunda expresión muestra que $\varphi(r_0)$ tiende a un límite $\varphi(\infty)$ cuando $r_0$ crece indefinidamente. La primera expresión de $\varphi(r_0)$ nos proporciona

$$|\varphi(r_0)| \leq \iint_s |u(x)a(x)|\frac{x_i\delta_i}{r_0}$$

de la cual

$$\int_0^\infty |\varphi(r_0)|dr_0 \leq \iiint_\Pi |u(x)a(x)|\delta x$$

L que resulta es $\varphi(\infty) = 0$. En otras palabras

$$\iiint_\Pi \left[u(y)\frac{\partial a(y)}{\partial y_i} + \frac{\partial u(y)}{\partial y_i}a(y)\right]\delta y = 0 \tag{11}$$

una formulación mas general de la anterior

$$\iiint_{\Pi-\bar{\omega}} \left[u(y)\frac{\partial a(y)}{\partial y_i} + \frac{\partial u(y)}{\partial y_i}a(y)\right]\delta y = -\iint_\sigma u(y)a(y)\delta y \tag{12}$$

Elegimos un dominio $\bar{\omega}$ como una esfera de radio infinitamente pequeño y centro $x$, y tomamos de la ecuación (12) $a(y) = \frac{1}{4\pi}\frac{\partial(\frac{1}{r})}{\partial y_i}$ y añadimos las relaciones de los valores $1, 2, 3$ de $i$ para obtener una identidad importante

$$u(x) = \frac{1}{4\pi}\iiint \frac{\partial(\frac{1}{r})}{\partial y_i}\frac{\partial u}{\partial y_i}\delta y \tag{13}$$

Ahora tomamos $a(y) = \frac{y_i - x_i}{r^2}$ en (11), y añadimos las mismas relaciones como en la expresión anterior

$$2\iiint_\Pi \frac{y_i - x_i}{r^2}\frac{\partial u}{\partial y_i}u(y)\delta y = -\iiint_\Pi \frac{1}{r^2}u^2(y)\delta y$$

aplicando la desigualdad de Schwarz al lado izquierdo obtenemos una desigualdad útil

$$\iiint_\Pi \frac{1}{r^2}u^2(y)\delta y \leq 4\iiint_\Pi \frac{\partial u}{\partial y_i}\frac{\partial u}{\partial y_i}\delta y \tag{14}$$

## 4.7. Cuasiderivadas

Sea $u^*(x)$ una infinidad de funciones cuadrado sumables con primeras derivadas cuadrado sumables continuas sobre $\Pi$. Suponemos que las derivadas $\frac{\partial u^*}{\partial x_1}, \frac{\partial u^*}{\partial x_2}, \frac{\partial u^*}{\partial x_3}$ convergen débilmente a la media sobre $\Pi$ a las funciones $U_{,1}, U_{,2}, U_{,3}$. Sea $U(x)$ una función medible (en el sentido de Lebesgue) definida en casi todo el espacio por

$$U(x) = \frac{1}{4\pi} \iiint_\Pi \frac{\partial(\frac{1}{r})}{\partial y_i} U_{,i} \delta y$$

tenemos entonces

$$\iiint_{\bar{\omega}} [u^*(x) - U(x)]^2 \, \delta x = -\iiint_\Pi \iiint_\Pi K_{ij}(y, y') \left[\frac{\partial u^*}{\partial y_i} - U_{,i}(y)\right] \left[\frac{\partial u^*}{\partial y'_j} - U_{,j}(y')\right] \delta y \delta y' \tag{15}$$

Siendo $r'$ la distancia entre $x$ e $y'$, también obtenemos lo siguiente

$$K_{ij}(y, y') = \frac{1}{16\pi^2} \iiint_{\bar{\omega}} \frac{\partial(\frac{1}{r})}{\partial y_i} \frac{\partial(\frac{1}{r'})}{\partial y'_j} \delta x$$

La expresión anterior para $K$, permite demostrar de manera sencilla que la integral

$$\iiint_\Pi \iiint_\Pi K_{ij}(y, y') K_{ij}(y, y') \delta y \delta y'$$

es finita, entonces el lado derecho de (15) está definido. Tiende a cero por (7). A partir de (11) podemos realizar la siguiente integral

$$\iiint_\Pi \left[U(y) \frac{\partial a}{\partial y_i} + U_{,i} a(y)\right] \delta y = 0 \tag{16}$$

Podemos dar la siguiente definición:

*Definición de Cuasiderivadas:* Consideremos dos funciones cuadrado sumable $U(y)$ y $U_{,i}(y)$ sobre $\Pi$. Podemos decir que $U_{,i}(y)$ es una cuasiderivada de $U(y)$ con respecto a $y_i$ cuando (16) se satisface. Recordando que en (16) $a(y)$ es cualquier función cuadrado sumable con primera derivada continua cuadrado sumable sobre $\Pi$

*Lema.* Supongamos que tenemos una infinidad de funciones continuas $u^*(x)$ con primeras derivadas continuas. Supongamos también que la integral $\iiint_\Pi u^{*2}(x) \delta x$ está acotada y para cada derivada $\frac{\partial u^*(x)}{\partial x_i}$ posee un límite débil en la media $U_{,i}(x)$ sobre $\Pi$. Entonces $u^*(x)$ converge en la media a función $U(x)$ para la cual $U_{,i}(x)$ son cuasiderivadas.

Siguiendo nuestra definición de cuasiderivadas tenemos que definir la *cuasidivergencia* $\Theta(x)$ de un vector $U_i(x)$ con componentes cuadrado sumable sobre $\Pi$. Cuando esta existe, es una función cuadrado sumable con

$$\iiint_\Pi \left[\Theta(y) a(y) + U_i(y) \frac{\partial a}{\partial y_i}\right] \delta y = 0 \tag{17}$$

## 4.8. Aproximación de una función medible mediante una sucesión de funciones regulares

Sea $\epsilon > 0$. Elegimos una función continua y positiva $\lambda(s)$ [2] definida para $0 \leq s$ e identicamente 0 para $s \leq 1$ y teniendo derivadas de cualquier orden tal que

$$4\pi \int_0^1 \lambda(\sigma^2) \sigma^2 d\sigma = 1$$

---

[2] Para aclarar ideas seleccionamos $\lambda(s) = A e^{\frac{1}{s-1}}$, $A$ es una constante arbitraria, $0 < s < 1$

Si $U(x)$ es sumable en todos los dominios $\bar{\omega}$, sea

$$\bar{U}(x) = \frac{1}{\epsilon^3} \iiint_\Pi \lambda\left(\frac{r^2}{\epsilon^2}\right) U(y) \delta y \tag{18}$$

$r$ es la distancia entre $x$ e $y$.
$\bar{U}(x)$ posee derivadas de todos los ordenes.

$$\frac{\partial^{l+m+n} \bar{U}(x)}{\partial x_1^l \partial x_2^m \partial x_3^n} = \frac{1}{\epsilon^3} \iiint_\Pi \frac{\partial^{l+m+n} \lambda\left(\frac{r^2}{\epsilon^2}\right)}{\partial x_1^l \partial x_2^m \partial x_3^n} U(y) \delta y \tag{19}$$

Sí $U(x)$ es acotado sobre $\Pi$, entonces claramente tenemos

$$min\ U(x) \leq \bar{U}(x) \leq máx\ U(x) \tag{20}$$

Si $U(x)$ es cuadrado sumable sobre $\Pi$ la desigualdad (3) aplicada a (18) se obtiene

$$\iiint_\Pi \bar{U}(x)^2 \delta x \leq \iiint_\Pi U^2(x) \delta x \tag{21}$$

Lo mismo aplicado a (19) prueba que las derivadas parciales de $\bar{U}(x)$ son cuadrado sumables sobre $\Pi$.
Finalmente tenemos lo siguiente, sí $U(x)$ y $V(x)$ son cuadrado sumable sobre $\Pi$

$$\iiint_\Pi \bar{U}(x) V(x) \delta x = \iiint_\Pi U(x) \bar{V}(x) \delta x \tag{22}$$

Sí $V(x)$ es continua $\bar{V}(x)$ tiende uniformemente a $V(x)$ sobre todo el dominio $\bar{\omega}$ cuando $\epsilon$ tiende a cero. *Lema:* Sea $U(x)$ cuadrado sumable sobre $\Pi$. $\bar{U}(x)$ converge fuertemente en la media a $U(x)$ sobre $\Pi$ cuando $\epsilon$ tiende a cero.

### 4.9. Algunos lemas de cuasiderivadas

Sea $U(x)$ cuadrado sumable en $\Pi$. Suponemos que para todas las funciones cuadrado sumables $a(x)$ poseen derivadas cuadrado sumables de todos los ordenes

$$\iiint_\Pi U(x) a(x) \delta x = 0$$

entonces

$$\iiint_\Pi U(x) \bar{U}(x) \delta x = 0$$

del cual uno puede tender a 0 $\epsilon$ de modo que

$$\iiint_\Pi U^2(x) \delta x = 0$$

donde $U(x)$ es por lo tanto 0 en cualquier parte.

*Lema 1:* Sa $U(x)$ un función que posee cuasiderivadas $U_{,i}(x)$. Entonces decimos que $\overline{\frac{\partial U(x)}{\partial x_i}} = \overline{U_{,i}(x)}$. Es suficiente probar que

$$\iiint_\Pi \overline{\frac{\partial U(x)}{\partial x_i}} a(x) \delta x = \iiint_\Pi U_{,i}(x) a(x) \delta x$$

Debido a que se puede deducir fácilmente de la expresión (18) que

$$\frac{\partial \overline{a(x)}}{\partial x_i} = \overline{\frac{\partial a(x)}{\partial x_i}}$$

# 5. Movimiento Infinitamente Lento

Las ecuaciones "linealizadas" de Navier-Stokes son las siguientes

$$\nu \Delta u_i(x,t) - \frac{\partial u_i(x,t)}{\partial t} - \frac{1}{\rho}\frac{\partial p_i(x,t)}{\partial x_i} = -X_i(x,t) \left[\Delta = \frac{\partial^2}{\partial x_k \partial x_k}\right] \qquad (23)$$

$$\frac{\partial u_j(x,t)}{\partial x_j} = 0$$

$\nu$ y $\rho$ son constantes dadas, $X_i(x,t)$ es un vector que representa todas las fuerzas externas presentes en el fenómeno, $p(x,t)$ es la presión, y $u_i(x,t)$ es la rapidez de las moléculas del fluido.

El problema dado en la teoría de fluidos viscosos es el siguiente: Construir para $t > 0$ una solución para (23) la cual tiene condición inicial $u_i(x,0)$.

Recordemos algunas la solución de este problema y algunas propiedades. Escribimos

$$W(t) = \iiint_\Pi u_i(x,t) u_i(x,t) \delta x$$

$$J_m^2(t) = \iiint_\Pi \frac{\partial^m u_i(x,t)}{\partial x_k \partial x_l \dots} \frac{\partial^m u_i(x,t)}{\partial x_k \partial x_l \dots} \delta x$$

$V(t) = $ Máximo de $\sqrt{u_i(x,t)u_i(x,t)}$ en el tiempo $t$ $D_m(t) = $ Máximo de la función $\left|\frac{\partial^m u_i(x,t)}{\partial x_1^h \partial x_2^k \partial x_3^l}\right|$ en el tiempo $t (l + k + h = m)$.

Haremos las siguientes aseveraciones. Las funciones $u_i(x,t)$ y sus primeras derivadas son continuas, $\frac{\partial u_j(x,t)}{\partial x_j}$, las cantidades $V(0)$ y $W(0)$ son finitas, $|X_i(x,t) - X_i(y,t)| < r^{1/2} C(x,y,t)$, donde $C(x,y,t)$ es una función continua, y $\iiint_\Pi X_i(x,t) X_i(x,t) \delta x$ es una función continua de $t$, o es menor que una función continua de $t$.

De ahora en adelante $A$ y $A_m$ denotarán constantes de índice $m$ las cuales no se les especificará un valor númerico.

## 5.1. Primer caso: $X_i(x,t) = 0$

La teoría de transferencia de calor proporciona la siguiente solución [Öz85] para (23):

$$u_i'(x,t) = \frac{1}{(2\sqrt{\pi})^3} \iiint_\Pi \frac{e^{-\frac{r^2}{4\nu t}}}{(\nu t)^{3/2}} u_i(y,0) \delta y; \qquad p'(x,t) = 0 \qquad (24)$$

La integral $u_i'(x,t)$ es uniformemente continua en $t$ para $t > 0$, y de esta integral podemos obtener lo siguiente

$$V(t) < V(0) \qquad (25)$$

Si $J_1(0)$ es finita, la desigualdad (14) y la desigualdad de Schwarz (1) aplicada a (24) nos proporciona la segunda cota para $V(t)$:

$$V^2(t) < 4 J_1^2(0) \frac{1}{(2\sqrt{\pi})^3} \iiint_\Pi \frac{e^{-\frac{r^2}{2\nu t}}}{(\nu t)^3} r^2 \delta y$$

es lo mismo expresarlo

$$V(t) < \frac{A J_1(0)}{(\nu t)^{1/4}} \qquad (26)$$

la desigualdad (3) aplicada (24) prueba:

$$W(t) < W(0) \qquad (27)$$

la integral $u'_i(x,t)$ es fuertemente continua en $t$ inclusive en $t = 0$. La desigualdad (3) aplicada a

$$\frac{\partial u'_i(x,t)}{\partial x_k} = \frac{1}{(2\sqrt{\pi})^3} \iiint_\Pi \frac{\partial}{\partial x_k} \left[ \frac{e^{-\frac{r^2}{4\nu t}}}{(\nu t)^{3/2}} \right] u_i(y,0)\delta y$$

prueba que

$$J_1(t) < J_1(0); \tag{28}$$

las primeras derivadas $\frac{\partial u'_i(x,t)}{\partial x_k}$ son fuertemente continuas en $t$, incluyendo $t = 0$ si $J_1(0)$ es finita.

Por razones análogas las derivadas de todos los órdenes de $u'_i(x,t)$ son uniformemente y fuertemente continuas en $t > 0$ y más precisamente

$$D_m(t) < \frac{A_m \sqrt{W(0)}}{(\nu t)^{\frac{2m+3}{4}}} \tag{29}$$

$$J_m(t) < \frac{A_m \sqrt{W(0)}}{(\nu t)^{\frac{m}{2}}} \tag{30}$$

### 5.2. Segundo Caso: $u'_i(x,0) = 0$

La solución fundamental de Oseen [Ose10], $T_{ij}(x,t)$ [3], provee la siguiente solución para (23);

$$u''_i(x,t) = \int_0^t dt' \iiint_\Pi T_{ij}(x-y,t-t')X_j(x,t')\delta y \tag{31}$$

$$p''(x,t) = -\frac{\rho}{4\pi} \frac{\partial}{\partial x_j} \iiint_\Pi \frac{1}{r} X_j(y,t)\delta y$$

Tenemos

$$|T_{ij}(x-y,t-t')| < \frac{A}{|r^2 + \nu(t-t')|^{3/2}} \tag{32}$$

$$\left| \frac{\partial^m T_{ij}(x-y,t-t')}{\partial x_1^h \partial x_2^k \partial x_3^l} \right| < \frac{A_m}{|r^2 + \nu(t-t')|^{\frac{m+3}{2}}}$$

Retomamos las integrales de (1) y (2), las aplicamos a (32) obteniendo entonces

$$\frac{\partial u''_i(x,t)}{\partial x_k} = \int_0^t dt' \iiint_\Pi \frac{\partial T_{ij}(x-y,t-t')}{\partial x_k} X_j(y,t)\delta y \tag{33}$$

esto prueba que las primeras derivadas $\frac{\partial u''}{\partial x_k}$ son fuertemente continuas en $t \geq 0$ y por lo tanto

$$J_1(t) < A \int_0^t \frac{dt'}{\sqrt{\nu(t-t')}} \sqrt{\iiint_\Pi X_i(x,t')X_i(x,t')\delta x} \tag{34}$$

Con esto hecho, añadimos anteriormente una hipótesis que aseveraba que el máximo de $\sqrt{X_i(x,t')X_i(x,t)}$ en el tiempo $t$ es una función continua de $t$. Entonces no es difícil demostrar de (31) que $u''_i(x,t)$ y $\frac{\partial u''}{\partial x_k}$ son uniformemente continuas en $t \geq 0$, y más precisamente

$$D_1(t) = A \int_0^t \frac{dt'}{\sqrt{\nu(t-t')}} max \sqrt{X_i(x,t)X_i(x,t)} \tag{35}$$

---

[3] $T_{ij} = \delta_{ij} \frac{1}{2\nu} \frac{E(r,t^{(0)}-t)}{t^{(0)}-t} + \frac{\partial^2 \Phi}{\partial x_j \partial x_k}, \Phi = \frac{1}{r} \int_0^t E(\alpha, t^{(0)}-t)d\alpha, E(r,t^{(0)}-t) = \frac{e^g}{}$

Complementamos la desigualdad (35) como sigue

$$\frac{\partial u_i''(x,t)}{\partial x_k} - \frac{\partial u_i''(y,t)}{\partial y_k} = \int_0^t dt' \iiint_{\overline{\omega}} \frac{\partial T_{ij}(x-z,t-t')}{\partial x_k} X_j(z,t')\delta z - \int_0^t dt' \iiint_{\overline{\omega}} \frac{\partial T_{ij}(y-z,t-t')}{\partial y_k} X_j(z,t')\delta z$$
$$+ \int_0^t dt' \iiint_{\Pi-\overline{\omega}} \left[ \frac{\partial T_{ij}(x-z,t-t')}{\partial x_k} - \frac{\partial T_{ij}(y-z,t-t')}{\partial y_k} \right] X_j(z,t')\delta z$$

$\overline{\omega}$ es el dominio de puntos con distancia menor $2r$ de $x$ e $y$. Si aplicamos la formula de diferencias finitas al corchete

$$\left[ \frac{\partial T_{ij}(x-z,t-t')}{\partial x_k} - \frac{\partial T_{ij}(y-z,t-t')}{\partial y_k} \right]$$

Se puede provar facilmente que:

$$\left| \frac{\partial u_i''(x,t)}{\partial x_k} - \frac{\partial u_i''(y,t)}{\partial y_k} \right| < A r^{1/2} \int_0^t \frac{dt'}{[\nu(t-t')]^{3/4}} max \sqrt{X_i(x,t)X_i(x,t)} \tag{36}$$

### 5.3. Caso General

Para obtener las soluciones $u_i(x,t)$ de (23) correspondiente al estado inicial $u_i(x,0)$, es suficiente superponer las soluciones particulares anteriores, tomando

$$u_i(x,t) = u_i'(x,t) + u_i''(x,t); \quad p(x,t) = p''(x,t)$$

Complementaremos la información de los párrafos anteriores, estableciendo que $u_i(x,t)$ es fuertemente continua en $t$.
La continuidad fuerte es evidente en el caso donde $X_i(x,t)$ sea cero fuera del dominio $\overline{\omega}$. Cuando $x$ se mueve infinitamente lejos, $u_i''(x,t)$, $\frac{\partial u_i''(x,t)}{\partial x_k}$ y $p(x,t)$ se aproximan a cero como $(x_i x_i)^{-3/2}, (x_i x_i)^{-2}$ y $(x_i x_i)^{-1}$ respectivamente y es suficiente con integrar

$$\nu u_i \Delta u_i - \frac{1}{2}\frac{\partial}{\partial t}(u_i u_i) - \frac{1}{\rho} u_i \frac{\partial p}{\partial x_i} = -u_i X_i$$

para obtener la *relación de disipación de energía*

$$\nu \int_0^t J_i^2(t)dt' + \frac{1}{2}(W(t) - W(0)) = \int_0^t dt' \iiint_{\Pi} u_i(x,t') X_i(x,t') \delta x \tag{37}$$

de la cual se puede obtener la siguiente desigualdad

$$\frac{1}{2}W(t) \leq \frac{1}{2}W(0) + \int_0^t dt' \sqrt{W(t')} \sqrt{\iiint_{\Pi} X_i(x,t) X_i(x,t) \delta x}$$

$W(t)$ posee al menos una solución $\lambda(t)$

$$\frac{1}{2}\lambda(t) = \frac{1}{2}W(0) + \int_0^t dt' \sqrt{\lambda(t')} \sqrt{\iiint_{\Pi} X_i(x,t) X_i(x,t) \delta x}$$

es lo mismo expresar

$$\sqrt{W(t)} \leq \int_0^t \sqrt{\iiint_{\Pi} X_i(x,t) X_i(x,t) \delta x} dt' + \sqrt{W(0)} \tag{38}$$

La expresión (38) se le denomina *desigualdad de disipación de energía de primer orden*.
Cuando $X_i(x,t)$ no es cero fuera del dominio $\overline{\omega}$, podemos puede aproximar las funciones $X_i(x,t)$ como una sucesión de funciones $X_i^*(x,t)$ que sean cero fuera del dominio $\overline{\omega}^*$, siempre y cuando se cumplan las relaciones (37) y (38). La expresión(37) muestra que $W(t)$ es continua. Luego entonces $u_i(x,t)$ es fuertemente continua en $t \geq 0$.

# Referencias